\documentclass[a4paper,11pt]{article}

\usepackage{stdformat}
\usepackage{standalone}
\usepackage{float}
\usepackage[shortlabels]{enumitem}
\usepackage[font=footnotesize, labelfont=bf]{caption}
\usepackage{subcaption}
\usepackage[affil-it]{authblk}

\numberwithin{equation}{section}
\numberwithin{figure}{section}
\numberwithin{table}{section}

\DeclareSymbolFont{matha}{OML}{txmi}{m}{it} 
\DeclareMathSymbol{\vv}{\mathord}{matha}{118}

\renewcommand{\thanks}[1]{\footnotetext{#1}}
\title{{\LARGE Multistage stochastic optimization for drayage procurement in container logistics using stochastic dual dynamic programming}}
\author[1,2]{Georgios Vassos}
\author[2]{Richard Lusby}
\author[3]{Pierre Pinson}

\affil[1]{Transported by Maersk, A.P. Moller - Maersk, 1098~Copenhagen~K, Denmark}
\affil[2]{Department of Technology, Management and Economics, Technical University of Denmark, 2800~Kgs.~Lyngby, Denmark}
\affil[3]{Dyson School of Design Engineering, Imperial College London, SW7~2AZ~London, United Kingdom}

\date{}

\begin{document}
\maketitle

\begin{abstract}
    \noindent Truckload procurement plays a vital role in integrated container logistics, particularly under the uncertainties of container flow and market conditions. We formulate the operational volume allocation problem in drayage procurement as a multistage stochastic transportation problem and solve it using stochastic dual dynamic programming (SDDP). We employ a multivariate count time series approach from the literature to model cargo flow dynamics, relaxing independence assumptions and capturing complex correlations. Our numerical experiments demonstrate the scalability of SDDP and its effectiveness in approximating high-quality policies across realistic problem instances. Sensitivity analyses highlight the significant impact of inflow uncertainties on costs, while spot market variability has a comparatively minor effect. Additionally, we propose an alternative stopping rule for SDDP iterations, balancing computational efficiency and solution fidelity. 
\end{abstract}

\section{Introduction}


Integrated container logistics companies coordinate the movement of goods across air, ocean, and inland modes of transport, leveraging containerization to reduce costs and enable seamless intermodal logistics \citep{Dekker2017}. Shipments are consolidated to optimize container space utilization, while the inland segment ensures connectivity between critical facilities like ports, rail yards, and distribution centers \citep{Yang2021}. Drayage supports this connectivity by linking intermodal terminals with customer locations, enhancing operational efficiency and sustainability \citep{Heggen2019}. A structured procurement process secures trucking capacity for short-distance transport, enabling adaptability to fluctuating shipment volumes and supply chain complexities. 

Two essential components of the procurement process are carrier selection and capacity planning \citep{Acocella2023}. Carrier selection involves evaluating potential carriers to identify reliable, long-term partners who support the organization throughout the planning horizon. Capacity planning, in turn, secures adequate trucking capacity with these trusted carriers to meet projected cargo flows in daily operations. An operational volume allocation policy completes the procurement strategy, serving as a rule for assigning container drayage execution in alignment with the capacity plan. This approach effectively integrates strategic and operational objectives by optimizing the operational policy within the framework of strategically secured capacities.


Service contracts enable shippers to secure carrier capacity aligned with forecast volumes while retaining the flexibility to adjust utilization based on evolving operational needs \citep{Acocella2023}. These contracts typically define rates per kilometer and total capacity allocations—often expressed in twenty-foot equivalent units (TEUs) for containerized cargo—providing a structured approach to managing transportation costs and mitigating demand uncertainty. By reducing reliance on potentially costly spot market alternatives, service contracts offer predictable rates and contribute to enhanced operational efficiency.



An integral component of the procurement process is the operational volume allocation policy, which provides a systematic approach for assigning drayage tasks to capacitated carriers. Previous research has highlighted that a well-designed policy optimizes resource utilization and supports reliable service delivery in dynamic transportation networks \citep{Rajapakshe2014, Triki2014}, and reduces costs \citep{Tsai2011, Othmane2019}. During the planning phase, this policy can be structured across temporal scales, such as quarterly, monthly, or daily, corresponding to planning horizons of 4, 12, or 365 stages annually. At each stage, the shipper monitors cargo flows within the drayage network and allocates volumes among carriers, who transport containers between entry and exit locations for consumption or onward shipping.


The drayage system, encompassing temporary container storage and transportation to exit points, operates under uncertain container flows and fluctuating spot market costs. In this study, we optimize volume allocation policies for drayage procurement by balancing commitments between contracted and spot carriers while minimizing costs associated with storage, dispatch delays, and container transportation. The optimization considers a set of capacitated carriers and assumptions about uncertain parameters, where decisions in carrier selection and capacity planning shape distinct operational policies. By computing the expected total cost of operations over the planning horizon under the optimal policy, we evaluate the broader procurement strategy, ensuring allocation decisions adapt to realized uncertainties while adhering to constraints like storage limits and carrier capacities.



We advance the study of truckload procurement by optimizing operational policies for drayage under container flow uncertainties. Our approach formulates the problem as a multistage stochastic transportation problem (MSTP) and employs SDDP to approximate optimal volume allocation policies. Unlike prior work by \cite{schmiedel2024} that primarily focuses on tactical-level decision integration, we emphasize operational policy optimization, supported by a novel method for multivariate count time series modelling proposed by \cite{Fokianos2020}, which captures correlations and relaxes traditional independence assumptions in cargo flow modelling. Sensitivity analysis reveals that inflow variations significantly affect objective values, while spot market uncertainty has a minor impact within the broader procurement strategy. We also identify diminishing marginal improvements in policy quality as SDDP iterations increase, leading to the proposal of an alternative stopping rule that balances computational efficiency with practical applicability. Numerical experiments demonstrate the scalability and fidelity of SDDP, bridging theoretical advancements and real-world operational challenges in transportation planning.


The paper is organized as follows: Section~\ref{sec:litrev} reviews the related literature on truckload procurement and multistage stochastic optimization. Section~\ref{sec:mmodel} introduces the mathematical model and methodology, detailing the formulation of the problem as an MSTP and the application of SDDP. Section~\ref{sec:numan} presents numerical experiments, including sensitivity analyses and computational performance evaluations. Finally, Section~\ref{sec:disco} concludes with a discussion of the findings and their practical implications.

\section{Literature review}\label{sec:litrev}


Recent literature defines the drayage problem as a complex logistical challenge that involves optimizing container movements between terminals, customer sites, and depots to meet diverse operational requirements. Drayage research frequently examines the efficient handling of import and export orders, ensuring containers are effectively transported to and from intermodal points as part of broader supply chains \citep{Escudero2021,Heggen2019}. Addressing the repositioning of empty containers is another key focus, with studies emphasizing strategies to reduce idle time and storage needs by directly moving containers between customer sites or consolidating them at depots \citep{Escudero2021,Song2017}. Recent work also considers the benefits of flexible timing and location assignments, which can significantly lower costs and improve scheduling flexibility under variable demand \citep{Moghaddam2020,Escudero2021}. Another critical aspect is the operational policy that allows trucks to decouple from containers during loading or unloading, which enhances truck utilization and allows the fleet to manage multiple tasks simultaneously. This approach, especially under \textit{separation mode} (explain), is analyzed by \cite{Song2017} for its impact on increasing drayage efficiency and reducing total route times. Together, these studies illustrate that modern drayage research seeks to balance efficiency, flexibility, and cost-effectiveness in container logistics through innovative routing and resource management solutions.


Literature on the drayage problem primarily addresses operational constraints such as container repositioning, minimizing idle time, optimizing timing and location assignments, and policies enabling trucks to decouple from containers during loading and unloading, all of which aim to enhance efficiency, flexibility, and cost-effectiveness. However, an often-overlooked factor impacting operational efficiency is drayage procurement. As a subset of the truckload market, the drayage sector is highly fragmented \citep{Gorman2023} and localized around intermodal hubs, with a substantial portion of its capacity managed by small operators specializing in short-haul container movements.

On the carrier side, capacity planning in truckload transportation is often reactive, with carriers arranging resources only after demand emerges from shippers, rather than utilizing advanced analytics \citep{Gorman2023}. This reactive approach leads to substantial volatility in carrier availability, while reliance on the spot market further exposes shippers to significant price fluctuations \citep{Budak2017}. To mitigate these risks, shippers frequently use reverse auctions \citep{Caplice2007} to secure long-term contracts that specify fixed carrier rates, ensuring stable capacity and shielding against market variation \citep{Sheffi2004}. Shippers with regular, periodic cargo flows--such as integrated logistics providers--typically use the spot market solely as a buffer to manage volatility in demand volumes \citep{Gorman2023}. This strategic approach helps shippers balance cost stability with flexible capacity options in response to fluctuating demand.


Academic literature on truckload procurement has primarily focused on auctions and spot market negotiations, with relatively fewer studies examining how operational procurement policies influence broader business performance \citep{Gorman2023,Lafkihi2019}. In a systematic review, \citet{Acocella2023} observed that most research emphasizes securing future transportation capacity through long-term contracts, typically using reverse auctions. However, the complexity of coordinating transportation services to meet varied logistical needs highlights opportunities for advanced analytic techniques to optimize routing, timing, and resource utilization \citep{Lafkihi2019}. \citet{Boada2020} explore dynamic procurement with partial demand information, focusing on cost and capacity alignment. \citet{Kantari2021} examine a mix of contract-based and on-demand sourcing to optimize shipment reliability and truck utilization under fluctuating demand. \citet{Hu2016} propose a bi-objective model that balances transportation costs and transit times, while \citet{Kuyzu2017} address collaborative lane covering using column generation for efficient routing. \citet{Oner2021} apply game theory to form stable coalitions, enhancing resource optimization and cost efficiency. In combinatorial auction contexts, \citet{Remli2013} and \citet{Zhang2014, Zhang2015} introduce robust and stochastic optimization models to manage shipment volume uncertainties, with \citet{Remli2019} extending this to address carrier capacity fluctuations. Continuing this trajectory, \citet{boujemaa2022} apply a SDDP approach to optimize carrier selection and shipment assignment under demand uncertainty. 

In this study, we frame truckload procurement operations as a multistage stochastic transportation problem (MSTP), minimizing the expected total inventory and transportation costs over a set planning horizon under constraints shaped by strategic procurement parameters. Addressing critical gaps in truckload procurement, we highlight capacity planning as a crucial, yet often overlooked, complement to carrier selection, introducing a tool for capacity optimization that also facilitates carrier evaluation. We employ a multivariate time series approach to model uncertainties in a manner that realistically reflects cargo dynamics, capturing interdependencies across flows to represent the exogenous network behavior more accurately. Finally, we assess the performance of SDDP in a case study grounded in practical size application, offering insights into its scalability and applicability in industry contexts.



\section{Mathematical model}\label{sec:mmodel}


Let $\mathcal{I}$ denote the set of origin inventories where containers are held upon entering the system boundary, and $\mathcal{J}$ denote the set of destination inventories where containers remain until they are either utilized within the system or dispatched externally. We define $m = |\mathcal{I} \times \mathcal{J}|$ as the total number of transportation lanes. The planning horizon is specified by the index set $\mathcal{T} = \{1, \dots, \tau\}$, with $\tau$ as the stopping time. At each time period~$t$, the state of the system is represented by the random vector $S_{t} : \Omega \to \mathcal{S}$, while interventions are modelled by the random decision vector $A_{t} : \Omega \to \mathcal{A}$, specifying drayage moves for each lane $(i,j) \in \mathcal{I} \times \mathcal{J}$ within the system allocated across $n$ selected carriers using both in-contract and spot sourcing. Note that \( \mathcal{S} \subset \mathbb{R}^{|\mathcal{I}|+2|\mathcal{J}|} \) as we consider one state variable for the stock level at each entry hub but two state variables one for the stock and one for the shortage level at each exit hub. Additionally, \( \mathcal{A} \subset [0, \infty)^{(2n) \times m} \) because transcations with each carrier can be either through contract or spot market. 

The cost incurred at each time $t$ following a decision $A_{t}$ consists of two components: the inventory holding and shortage costs, $h_{t} : \mathcal{S} \to [0, \infty)$, and the drayage transportation cost, $f_{t} : \mathcal{A} \to [0, \infty)$. Additionally, the exogenous state captures the cargo flows entering and exiting the system, denoted $Q_{t} : \Omega \to \mathcal{Q} \subset [0, \infty) ^ {|\mathcal{I}|}$ and $D_{t} : \Omega \to \mathcal{D} \subset [0, \infty) ^ {|\mathcal{J}|}$, as well as spot rates $R_{t} : \Omega \to \mathcal{R} \subset [0, \infty)^{n \times m}$. Additional design parameters include stock limits for both origin and destination inventories \( \overline{S} = \{ \overline{S}_{i} : i \in \mathcal{I}, \overline{S}_{j} : j \in \mathcal{J}  \} \), shortage limits \( \underline{S} = \{ \underline{S}_{j} : j \in \mathcal{J} \} \) specific to destinations, and a strategic sourcing arrangement characterized by execution rates per transportation lane $w_{t} \in [0, \infty) ^ {n \times m}$ and maximum capacity allocations for each strategic carrier $x_{t} \in [0, \infty) ^ {2n}$.

For brevity, let \(\mathcal{L}\) denote the set of all lanes \(\mathcal{I} \times \mathcal{J}\). A carrier \(k \in \mathcal{K} = \{ 1, \dots, 2n \}\) can serve all lanes through spot procurement, when $n < k \le 2n$, at a variable rate, but it is limited to serving a subset \(\mathcal{L}(k) \subset \mathcal{L}\), when $1 \le k \le n$, at a fixed rate, as specified in an awarded contract. The state vector is expanded as $S_{t} = \{ S_{i,t} : i \in \mathcal{I}, (S_{j,t}^{+}, S_{j,t}^{-}) : j \in \mathcal{J} \}$, with $S_{j,t}^{+} = \max\{S_{j,t}, 0\}$ and $S_{j,t}^{-} = \max\{-S_{j,t}, 0\}$, and the decision vector is:
\[
    A_{t} = \big\{ \underbrace{(A_{i,j,t}^{k}:(i,j)\in\mathcal{L}(k),k = 1, \dots, n)}_{\text{contract}}, \underbrace{(A_{i,j,t}^{k} :(i,j)\in\mathcal{L},k = n + 1, \dots, 2n)}_{\text{spot}} \big\}
\]
Similarly, the capacity vector is given by:
\[
    x_{t} = \big\{ \underbrace{(x_{t}^{k} : k = 1, \dots, n)}_{\text{contract}}, \underbrace{(x_{t}^{k} : k = n+1, \dots, 2n)}_{\text{spot}} \big\}
\]
We assume a linear functional form for \( h_{t} \), with time homogeneous coefficients \((\gamma_{i} : i \in \mathcal{I})\) and \((\alpha_{j} : j \in \mathcal{J})\) for the costs of holding container stocks at entry and exit locations, respectively, and \((\beta_{j} : j \in \mathcal{J})\) for the cost of unmet outflow at exit locations.


Following the approach of \cite{boujemaa2022}, we formulate the multistage transportation optimization problem with a linear objective function assuming $w_{t} = w$ for all $t \in \mathcal{T}$, structured around specified contract rates and carrier capacities, as follows:
\begin{align}
     \underset{S_{t}, A_{t}\, :\, t \in \mathcal{T}}{\text{minimize}}&\ \sum_{t\in\mathcal{T}}\sum_{i\in\mathcal{I}}\gamma_{i}S_{i,t}+\sum_{t\in\mathcal{T}}\sum_{j\in\mathcal{J}}\left(\alpha_{j}S_{j,t}^{+}+\beta_{j}S_{j,t}^{-}\right)\nonumber\\
     &\hspace{15pt}+\sum_{t\in\mathcal{T}}\sum_{k = 1}^{n}\sum_{(i,j)\in\mathcal{L}(k)} w_{i,j}^{k} A_{i,j,t}^{k}+\sum_{t\in\mathcal{T}}\sum_{k = 1}^{n}\sum_{(i,j) \in \mathcal{L}} R_{i,j,t}^{k} A_{i,j,t}^{n+k}\\
    \text{subject to:}\hspace{-25pt} & \nonumber \\
     & \sum_{(i,j)\in\mathcal{L}(k)}A_{i,j,t}^{k}\le x_{t}^{k} && \forall\,t\in\mathcal{T},k\in\mathcal{K}\label{lp:cf}\\
     & \sum_{k \in \mathcal{K}}\sum_{j\in\mathcal{J}}A_{i,j,t}^{k} \le S_{i,t} + Q_{i,t} && \forall\,t\in\mathcal{T},i\in\mathcal{I}\label{lp:storlim1}\\
     & \sum_{k \in \mathcal{K}}\sum_{i\in\mathcal{I}}A_{i,j,t}^{k} +  S_{j,t}^{+} \le \overline{S}_{j} && \forall\,t\in\mathcal{T},j\in\mathcal{J}\label{lp:storlim2}\\
     & \sum_{k \in \mathcal{K}}\sum_{(i,j)\in\mathcal{L}(k)}A_{i,j,t}^{k}-S_{i,t}+S_{i,t+1}=Q_{i,t} && \forall\,t\in\mathcal{T},i\in\mathcal{I} \label{eq:trans1} \\
     & \sum_{k \in \mathcal{K}}\sum_{(i,j)\in\mathcal{L}(k)}A_{i,j,t}^{k}+S_{j,t}^{+}-S_{j,t}^{-}-S_{j,t+1}^{+}+S_{j,t+1}^{-}=D_{j,t} && \forall\,t\in\mathcal{T},j\in\mathcal{J} \label{eq:trans2} \\
     & S_{i,t} \ge 0 && \forall\,t\in\mathcal{T},i\in\mathcal{I}\\
     & S_{j,t}^{+},S_{j,t}^{-} \ge 0 && \forall\,t\in\mathcal{T},j\in\mathcal{J}\label{lp:cx}\\
     & A_{i,j,t}^{k} \ge 0 && \hspace{-60pt} \forall\,t \in \mathcal{T},k \in \mathcal{K}, (i,j)\in\mathcal{L}(k)\label{lp:cxx}
\end{align}
At a given time $t$, we subsume all uncertainty under the general parameter $\xi_{t} = (R_{t}, Q_{t}, D_{t}) : \Omega \to \Xi$. We let $\Xi = \mathcal{R} \times \mathcal{Q} \times \mathcal{D}$ be the support of the uncertain parameter and $\Xi^{\tau} = \Xi \times \cdots \times \Xi$ the scenario space. The topology of the scenario space induced by the probability measure $P$, of the ambient probability space $(\Omega, \mathcal{F}, P)$, is what we need to control in order to make our solution tractable. In general, the nature of such multivariate panel data poses significant challenges.

The multistage linear program outlined above cannot be solved directly due to the presence of random variables. Stochastic programming addresses this challenge by optimizing under uncertainty, using the probability distribution of these random variables to inform decision-making. The objective function is a functional of the distribution that generates an observation:
\begin{equation}\label{eq:obs}
    O_{\tau} = (S_{1},A_{1},S_{2},A_{2},\dots,S_{\tau},A_{\tau},S_{\tau+1})
\end{equation}
For a given scenario $\xi \doteq \xi_{[\tau]} = (\xi_{1},\dots,\xi_{\tau}) \in \Xi^{\tau}$, the probability density of $O_{\tau}$ is commonly factorized as the product:
\begin{equation}
    p(o_{\tau} \mid \xi) = p_{1}(s_{1})\prod_{t=1}^{\tau}p_{t}(s_{t+1} \mid s_{t}, a_{t}, \xi_{[t]})\,\pi_{t}(a_{t} \mid s_{t}, \xi_{[t]})
\end{equation}
where \(\xi_{[t]} = (\xi_{1}, \dots, \xi_{t})\) represents the exogenous state variable up to time \(t\). We identify three key parameter categories: the initial state density \(p_{1}(s_{1})\), the transition dynamics \(p_{t}(s_{t+1} \mid s_{t}, a_{t}, \xi_{[t]})\) for each \(t \in \mathcal{T}\), and the decision policy \(\pi_{t}(a_{t} \mid s_{t}, \xi_{[t]})\) across all \(t \in \mathcal{T}\). Typically, the initial state \(s_{1}\) is predetermined within \(\mathcal{S}\), and the policy emerges from the solution as the optimal decision vector. Transition dynamics are governed by constraints \eqref{eq:trans1} and \eqref{eq:trans2}. By extending the state vector \(S_{t}\) with the exogenous state \(\xi_{[t]}\), we encapsulate all stochastic influences within \(\xi_{t}\), allowing us to reformulate the dynamics as:
\begin{equation}
    p_{t}(s_{t+1} \mid s_{t}, a_{t}, \xi_{[t]}) = \delta(s_{t+1} \mid s_{t}, a_{t}, \xi_{t})
\end{equation}
Dirac's delta notation is mildly abused to highlight that $s_{t+1}$ is known with certainty given knowledge of $s_{t}$, $a_{t}$, and $\xi_{t}$. Thereby, the only effective nuisance parameter is the distribution of a random scenario $\xi$. 

Allowing $h_{t}$ to be time homogeneous by letting $h_{t} = h$ for all $t \in \mathcal{T}$, we define the stage specific inventory cost function as:
\begin{equation}
    h(s_{t}) = \sum_{i\in\mathcal{I}}\gamma_{i}s_{i,t}+\sum_{j\in\mathcal{J}}\left(\alpha_{j}s_{j,t}^{+}+\beta_{j}s_{j,t}^{-}\right)
\end{equation}
Moreover, we define the transportation cost function $f_{t} : \mathcal{S} \times \mathcal{A} \to [0, \infty)$ by:
\begin{equation}
    f_{t}(a_{t}) = \sum_{k = 1}^{n}\sum_{(i,j)\in\mathcal{L}(k)} w_{i,j}^{k} a_{i,j,t}^{k}+\sum_{k = 1}^{n}\sum_{(i,j) \in \mathcal{L}} r_{i,j,t}^{k} a_{i,j,t}^{n+k}
\end{equation}
We define the objective function of the multistage stochastic linear program to be the expected total inventory and transportation cost over the finite planning horizon:
\begin{equation}
    E_{O_{\tau}}\left[ \sum_{t \in \mathcal{T}} \big( h(S_{t}) + f_{t}(A_{t}) \big) \right] = E_{\xi}\left[ \sum_{t \in \mathcal{T}} \big( h(S_{t}) + f_{t}(A_{t}) \big) \right]
\end{equation}
The optimization is constrained by equations and inequalities \eqref{lp:cf}-\eqref{lp:cxx}. For any $\xi_{t} \in \Xi$, we let $\mathcal{O}_{t}(\xi_{t}) \subset \mathcal{S} \times \mathcal{A} \times \mathcal{S}$ denote the feasible set of $(s_{t}, a_{t}, s_{t+1})$, at each $t \in \mathcal{T}$.

For a given scenario $\xi \in \Xi^{\tau}$, we observe the instance that solves the linear program:
\begin{equation}
    o_{\tau}(\xi) = \argmin \left\{ E_{O_{\tau} \mid \xi}\left[ \sum_{t \in \mathcal{T}} \big( h(S_{t}) + f_{t}(A_{t}) \big) \right] : o_{\tau}\in\bigcap_{t\in\mathcal{T}}\mathcal{O}_{t}(\xi_{t}) \right\}
\end{equation}

The idea with stochastic optimization is to obtain the realization that is optimizing over a functional of the uncertainty distribution. For example, we may consider the expected total cost across the entirety of the scenario space:
\begin{equation}
    o_{\tau} = \argmin\left\{ E_{\xi}\left[ \sum_{t \in \mathcal{T}} \big( h(S_{t}) + f_{t}(A_{t}) \big) \right] : o_{\tau} \in \bigcap_{\xi \in \Xi^{\tau}} \bigcap_{t\in\mathcal{T}} \mathcal{O}_{t}(\xi_{t})\right\}
\end{equation}
Perhaps we are interested in minimizing the expected total cost for a subset of the scenario space $\Xi_{N} \subset \Xi^{\tau}$, where $|\Xi_{N}| = N$, which corresponds to specific uncertain conditions.

In multistage stochastic programming, an optimal solution can be obtained through backward induction, a category of dynamic programming. Here, the principle of optimality is applied by defining the terminal state cost and iterating backward from \( t = \tau \) to \( t = 1 \), where at each stage \( t \) the expected value is maximized over all feasible decisions. Specifically, the immediate reward function at each stage is given by:
\begin{align}
     u_{t}(S_{t}, A_{t}) &= - h(S_{t}) - \min_{A_{t}}\ f_{t}(A_{t})\\
     &\hspace{-25pt}\text{subject to: } \eqref{lp:cf}, \eqref{lp:storlim1}, \eqref{lp:storlim2},\text{ and }\eqref{lp:cxx} \nonumber
\end{align}
This random reward function \( u_{t}(S_{t}, A_{t}) \doteq u_{t}(s_{t}, a_{t};\xi_{t}) \) depends on scenario \( \xi_{t} \), which encapsulates uncertainty in state and decision variables at each stage.


For each $s_{t} \in \mathcal{S}$, we write the well-known backward induction formula:
\begin{equation}
    V_{t}(s_{t}) = \max_{a_{t} \in \mathcal{A}_{t}} \left\{ E_{\xi_{t}} u_{t}(s_{t}, a_{t}) + E_{S_{t+1}} V_{t+1}(s_{t},a_{t}) \right\}
\end{equation}
where \( \mathcal{A}_{t} = \mathcal{A}(s_{t})\) is a subset of \( \mathcal{A} \), \( E_{S_{t+1}} V_{t+1} : \mathcal{S} \times \mathcal{A} \to \mathbb{R} \) and
\begin{equation}
\begin{aligned}
    E_{S_{t+1}} V_{t+1}(s_{t}, a_{t}) &= E_{S_{t+1}} [V(S_{t+1})]\\
    &= \int_{\mathcal{S}} V_{t+1}(s_{t+1}) \, p_{t}(ds_{t+1} \mid s_{t}, a_{t}) \\
    &= \int_{\Xi} \int_{\mathcal{S}} V_{t+1}(s_{t+1}) \, \delta(ds_{t+1} \mid s_{t}, a_{t}, \xi_{t}) \, p_{t}(d\xi_{t}) \\
    &\doteq E_{\xi_{t}} V_{t+1}(S_{t+1})
\end{aligned}
\end{equation}
Hence, the backward induction is finally given by:
\begin{equation}\label{eq:optim}
    V_{t}(s_{t}) = \max_{a_{t} \in \mathcal{A}_{t}} \left\{ E_{\xi_{t}} \big[ u_{t}(s_{t}, a_{t}) + V_{t+1}(S_{t+1}) \big] \right\}
\end{equation}
Computing \( V_{t}(s_{t}) \) based on \eqref{eq:optim} for each \( t \in \mathcal{T} \) and \( s_{t} \in \mathcal{S} \), with \( |\mathcal{A}| \) evaluations of the expectation over \( \Xi \), faces dimensionality challenges. The state space \( \mathcal{S} \), action space \( \mathcal{A} \), and uncertainty support \( \Xi \) each grow exponentially with the number of inventories, while \( \mathcal{A} \) and \( \Xi \) further depend on the number of carriers. To alleviate the dimensionality burden in \( \Xi \), subsampling \( N \) scenarios \( \Xi_{N} \) from \( \Xi^{\tau} \) with \( N \ll |\Xi^{\tau}| \) is an effective approach. However, more advanced techniques are required to handle the dimensionality of the state and decision spaces.

\subsection{Stochastic Dual Dynamic Programming}

SDDP is a state-of-the-art algorithm for solving multistage stochastic programming problems, especially in situations where decisions must be made sequentially under uncertainty. Developed initially for applications like energy planning, SDDP is a sampling-based extension of the nested decomposition method that focuses on managing high-dimensional decision processes \citep{Dowson2021}. It decomposes large-scale stochastic problems into smaller subproblems using value-function approximations, making them computationally feasible to solve.

In SDDP, stages represent decision points where uncertainty is resolved progressively. At each stage, SDDP computes decision rules by approximating the cost-to-go function, iterating between forward passes (sampling possible outcomes) and backward passes (updating approximations). A thorough analysis of the SDDP algorithm and its convergence properties is provided by \citet{Shapiro2011}. This approach enables efficient handling of scenarios with a vast number of potential future states.

SDDP formulates the problem as a sequence of two-stage stochastic programs. At each stage $t$, the algorithm first calculates the reward based on the current decision and observed realization of uncertainty, capturing the immediate cost. Then, it estimates the expected future cost by incorporating decisions that account for scenario-based outcomes at the next stage $t+1$. This future cost is constrained by the first stage decision and approximated iteratively using backward and forward passes across the sampled scenarios \citep{Shapiro2011}. The immediate reward function is now defined as follows:
\begin{equation}
    c_{t}(s_{t},a_{t};\xi_{t}) = h(s_{t}) + f_{t}(a_{t};\xi_{t})
\end{equation}
According to \citet{Dowson2021}, in SDDP, the principle of optimality is expressed, for all \(s_{t} \in \mathcal{S}\) and \( \xi_{t} \in \Xi \), as:
\begin{equation}
    J_{t}(s_{t},\xi_{t}) = \min \left\{ c_{t}(s_{t},a_{t};\xi_{t}) + E_{\xi_{t+1}} [J_{t+1}(s_{t+1},\xi_{t+1})] : (s_{t}, a_{t}, s_{t+1}) \in \mathcal{O}_{t}(\xi_{t}) \right\}
\end{equation}
The algorithm based on the above equation is intractable because solving it requires evaluating the expectation over all possible next-stage scenarios at each iteration, which leads to an exponential growth in computational complexity that cannot be managed directly. SDDP addresses the intractability by approximating the expectation and iteratively refining this approximation with cuts, which approximate the cost function from below \citep{Dowson2021}. This approach reduces computational complexity by avoiding the need to evaluate all scenarios at each iteration.

\subsection{General multivariate uncertainty model}\label{sec:nuisance}

In this section, we introduce a robust approach for modelling multidimensional uncertainty. This methodology applies to panel data structures, where inflows occur across designated entry points and outflows across corresponding exit points. The approach effectively captures both the autocorrelation within marginal flows and the cross-correlation across panels, facilitating an integrated analysis of dependencies across multiple dimensions.


We define marginal flows as nonstationary autoregressive interval counts $(Q_{i,t} : i \in \mathcal{I}, D_{j,t} : j \in \mathcal{J}) \in \mathbb{N}^{|\mathcal{I}|+|\mathcal{J}|}$, for all $t \in \mathcal{T}$, with intensities specified by:
\begin{align}
    \lambda_{i}(t) &= \kappa_{i,0}+\kappa_{i,1}\lambda_{i}(t-1) + \kappa_{i,2}Q_{i,t-1} + \kappa_{i,3}{}^{\top}Z_{i,t} && \forall\, i \in \mathcal{I} \\
    \lambda_{j}(t) &= \kappa_{j,0}+\kappa_{j,1}\lambda_{j}(t-1) + \kappa_{j,2}D_{j,t-1} + \kappa_{j,3}{}^{\top}Z_{j,t} && \forall\, j \in \mathcal{J} 
\end{align}
where $\{ Z_{i,t} : i \in \mathcal{I} \}$ and $\{ Z_{j,t} : j \in \mathcal{J} \}$ represent exogenous factors, such as seasonal trends. At each time stage $t$, inflows can be simulated by drawing samples from $\text{Poisson}(\lambda_{i}(t))$ for all $i \in \mathcal{I}$, and outflows by drawing samples from $\text{Poisson}(\lambda_{j}(t))$ for all $j \in \mathcal{J}$.

Our methodology for constructing a joint distribution from marginal distributions with embedded correlations leverages copula modelling techniques \citep{Hofert2019}. Following the algorithm outlined in \citep{Fokianos2020}, we generate multivariate Poisson counts by sampling correlated uniform random variables from a copula, transforming these to exponential variables, and counting terms until reaching a threshold value of 1. This approach preserves the Poisson marginal distributions while introducing correlations across dimensions through the copula framework.

To capture these correlations, we construct a Gaussian copula by first sampling from a multivariate normal distribution \( \mathcal{N}(0_{|\mathcal{I}| + |\mathcal{J}|}, \Sigma) \), where \( 0_{|\mathcal{I}| + |\mathcal{J}|} \) is a vector of \( |\mathcal{I}| + |\mathcal{J}| \) zeroes and \( \Sigma \in \mathbb{R}^{(|\mathcal{I}| + |\mathcal{J}|) \times (|\mathcal{I}| + |\mathcal{J}|)} \) is the covariance matrix encoding the desired correlations. We then transform each component of this sample to a uniform random vector by applying the marginal cumulative distribution function of \( \mathcal{N}(0, \sigma_{ii}^2) \) for variables corresponding to inflows and \( \mathcal{N}(0, \sigma_{jj}^2) \) for outflows. This process yields a set of correlated uniform random variables, forming the basis for generating multivariate Poisson counts.


\section{Numerical analysis}\label{sec:numan}

The computational experiments were conducted on a MacBook Pro equipped with an Apple M1 Pro chip, featuring a 10-core CPU (comprising 8 performance cores and 2 efficiency cores) with a clock speed of 3.2 GHz. The system had 32 GB of unified memory and was running macOS Sonoma Version 14.7.1. All computations were executed using Julia version 1.11.1. The stochastic dual dynamic programming models were implemented using the \texttt{SDDP.jl} package (version 1.9.0) and solved with the Gurobi Optimizer via the \texttt{Gurobi.jl} interface (\texttt{Gurobi.jl} package version 1.4.0). The specific version of the Gurobi Optimizer used was 11.0.2 build v11.0.2rc0 (\texttt{mac64[arm]} - Darwin 23.6.0 23H222).

\subsection{Optimal volume allocation policy}



We illustrate the optimal operational policy using a simple problem instance involving two entry and two exit hubs, each with a maximum capacity of 100 TEUs. Two carriers can fulfill transportation demands through contracts or spot market bids across all four lanes. Storage costs are \$20/TEU at entry hubs and \$10/TEU at exit hubs, while unmet outflow incurs a shortage cost of \$30/TEU. Transportation costs for winning bids are uniformly distributed between \$7.0 and \$9.9/TEU, and spot rates range between \$3.5 and \$8.0/TEU.

\begin{figure}
    \centering
    \includegraphics[width = \textwidth]{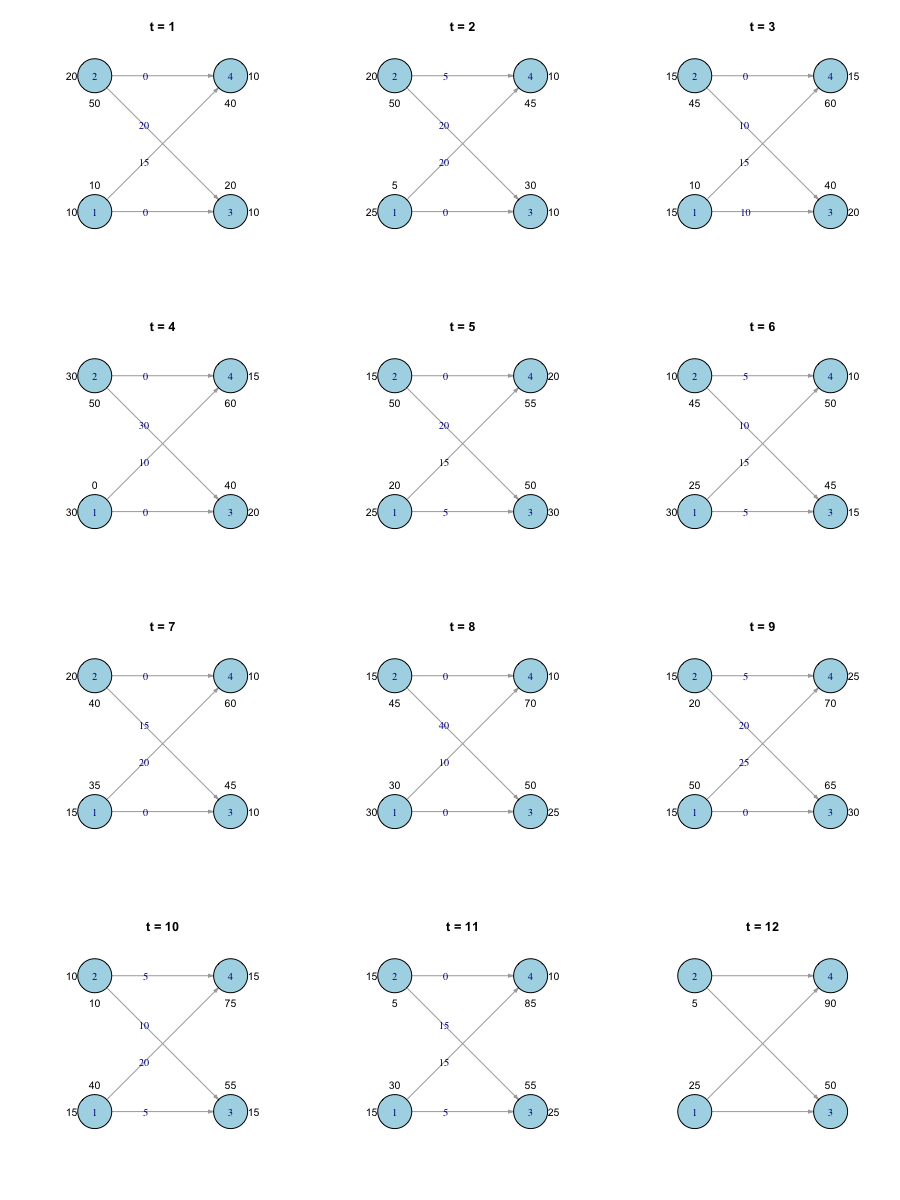}
    \caption{Optimal drayage operations in a system with 2 entry hubs and 2 exit hubs over 12 time stages, illustrating stock levels, inflows, outflows, and volume allocations between hubs.}
    \label{fig:2x2opt}
\end{figure}

Inflows and outflows over the planning horizon follow independent discrete uniform distributions over \(\{10, 15, 20, 25, 30\}\), and contracted carrier capacity per period is drawn uniformly from \(\{10, 15, 20\}\). Figure~\ref{fig:2x2opt} illustrates the optimal volume allocation for a deterministic instance, where nodes 1 and 2 represent entry hubs and nodes 3 and 4 represent exit hubs. Stock levels are shown above and below the nodes, inflows are displayed to the left of entry hubs, and outflows to the right of exit hubs. Arrows indicate lanes, with the numbers representing total volumes shipped between hubs at each time stage.

\begin{figure}
    \centering
    \includegraphics[width = 0.7\textwidth]{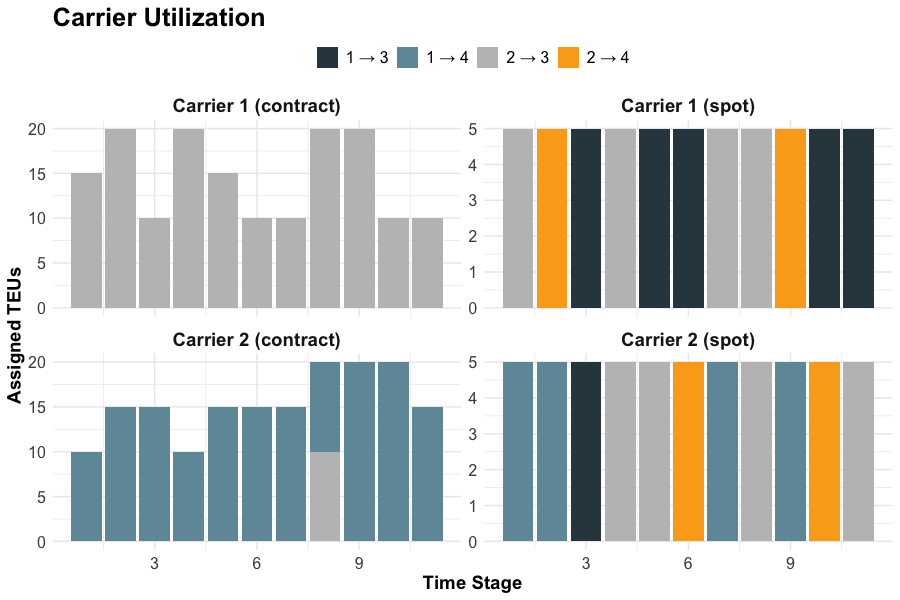}
    \caption{Optimal volume allocation between carriers and lanes.}
    \label{fig:2x2util}
\end{figure}

Figure~\ref{fig:2x2util} details carrier utilization for contract and spot segments at each time stage, breaking down the volumes in Figure~\ref{fig:2x2opt} into TEU-specific assignment decisions. The lane from location 2 to 3 was predominantly served by Carrier 1, while the lane from 1 to 4 was exclusively served by Carrier 2. This pattern reflects the cost structure, where each winning bid establishes a contract with a fixed transportation cost per TEU over the planning horizon. Additionally, factors such as carrier fleet locations influence the fixed transportation costs per lane. As a result, Carrier 1 offered more favorable contract terms for the lane from location 2 to 3, leading to minimal utilization of Carrier 2 on this lane during time stage 8, after exhausting Carrier 1's contracted capacity and all spot capacity. Overall, the spot market was leveraged opportunistically, exploiting cost variations to buffer operational costs.

For this specific problem instance, we observe full carrier utilization to maintain minimum stock levels at entry hubs, where higher storage costs make this strategy essential. While contracted carrier capacity is relatively more expensive than spot capacity due to the cost of securing guaranteed volumes, a small buffer remains available for leveraging lower-cost opportunities in the spot market. The challenges associated with securing carrier capacity over the planning horizon are discussed in \citep{Acocella2023}. Shippers are likely to accept higher contract rates to secure carrier capacity, while spot market can be utilized selectively to enhance operational efficiency.

The small-scale instance effectively illustrated the execution of the volume allocation policy; however, real-world challenges in container logistics are significantly more complex. In practice, container logistics companies often manage scenarios with a small number of high-capacity entry and exit locations and a large number of low-capacity carriers, reflecting the fragmented nature of the drayage business. The majority of trucking capacity is secured through contractual agreements to ensure frequent and continuous periodic moves, while a small spot market buffer is typically used to boost operational efficiency.

\subsection{Practical instance specifications}\label{sec:realinst}


The construction of problem instances that match the scale of real-world drayage systems involves defining a 12-period drayage system with parameters that emulate real-world complexities. The system comprises 6 entry hubs and 6 exit hubs, alongside 20 carriers competing for 1–2 bids each from a pool of 10 bids, where each bid includes 6–18 origin-destination (OD) pairs. Initial TEU stocks are sampled uniformly between 0 and 500 at entry hubs and between 0 and 1,000 at exit hubs, ensuring no initial shortages. Entry and exit capacities are fixed at 10,000 TEUs to reflect operational constraints, while carrier capacities range between 400 and 800 TEUs per period for contracts, supplemented by spot buffers fixed at 40 TEUs per period.

Cost parameters capture essential trade-offs, including storage costs (\$20.0/TEU at entry and \$10.0/TEU at exit hubs), shortage penalties (\$30.0/TEU), contract rates (\$6.0--\$8.0/TEU), and spot market rates (\$3.0--\$9.0/TEU). Demand and supply dynamics are simulated by drawing inflows and outflows at each hub uniformly at random between 1,000 and 3,000 TEUs per period. This configuration ensures the inclusion of bid allocations, dynamic flows, and cost structures under capacity constraints, enabling the creation of realistic and representative problem instances.

\subsubsection{Sensitivity analysis}

In their recent work, \cite{schmiedel2024} and \cite{boujemaa2022} pointed out that much of the existing literature on multistage transportation often overlooks the stochastic nature of the problem. While they addressed this gap by considering uncertainty in demand (outflow), their approach relied on assumptions of independence across time stages, products, and distribution centers (exit hubs). In this study, we extend this line of inquiry by investigating additional sources of uncertainty--namely, inflow (supply) and spot rates. To this end, we perform a sensitivity analysis on the objective value (total cost) using 100 problem instances generated as outlined in Section~\ref{sec:realinst}, aiming to assess whether these factors warrant consideration when modelling the problem.

\begin{figure}
    \centering
    \includegraphics[width = 0.95\textwidth]{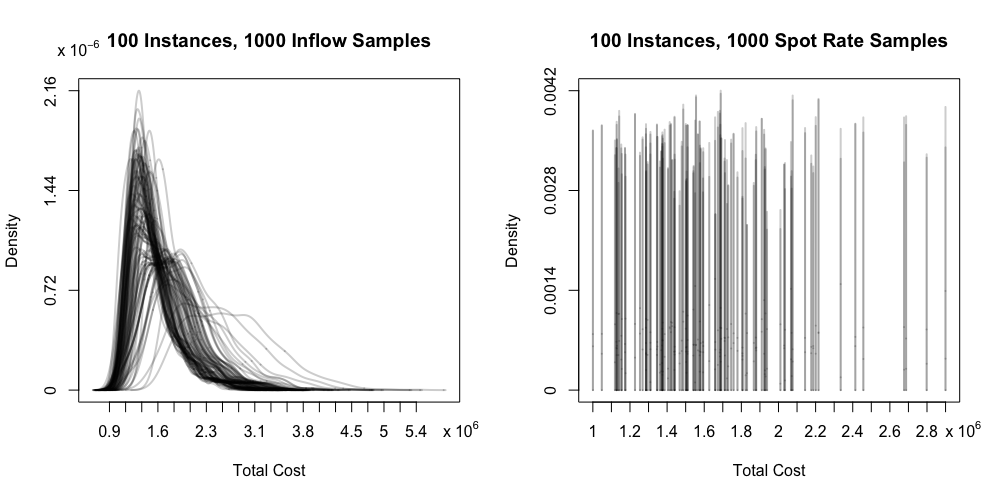}
    \caption{Distribution of the objective function across 1000 examples of inflow and spot rates for 100 instances of the 12-stage drayage problem, featuring 6 entry hubs, 6 exit hubs, and 20 carriers operating in both contract and spot markets.}
    \label{fig:sensv}
\end{figure}

For each instance, the inflow vector--comprising 72 elements with values uniformly distributed between 1,000 and 3,000--is sampled 1,000 times. Figure~\ref{fig:sensv} demonstrates the impact of inflow variations on the total cost across the 100 instances. The distinct density curves reveal substantial variability in the objective value due to changes in the inflow vector, underscoring the importance of incorporating inflow uncertainty into the analysis.

The right panel of Figure~\ref{fig:sensv} explores the sensitivity of total cost to spot rate variability across the same 100 instances. For each instance, 1,000 spot rate vectors are sampled uniformly between \$3.0/TEU and \$9.0/TEU. Despite significant spot rate fluctuations, the total cost exhibits minimal variability within individual instances. This suggests that reliance on the spot market remains limited under the evaluated procurement policies. Consequently, spot rate variability has a negligible impact on the total cost, shifting analytical attention to more influential stochastic factors, such as inflow uncertainty.

\subsection{Computational performance of the SDDP}

The computational performance of the SDDP algorithm in solving the multistage stochastic transportation problem is evaluated across 100 problem instances generated as outlined in Section~\ref{sec:realinst}. Metrics include bias, defined as the expected difference between the upper and lower bounds of the value function approximation during the forward and backward passes, following \citep{schmiedel2024}; the confidence interval (CI) ratio, representing the relative width of the CI for the expected objective value during training, normalized by its mean; and the average execution time required to solve a single problem instance. The operational volume allocation policy is assessed via regret analysis, measuring the relative deviation of its expected objective value from the deterministic optimum using out-of-bag scenarios from the same distribution.

The distribution of uncertain parameters for training the operational policy using SDDP and testing it on out-of-bag examples is modelled as a copula-based distribution with Poisson marginals. A Gaussian copula is employed to introduce correlation between the marginals, with a constant rate parameter \(\lambda_{i}(t) = \lambda_{j}(t) = 2,000\) for all \(i \in \mathcal{I}\), \(j \in \mathcal{J}\), and \(t \in \mathcal{T}\). This approach relaxes the common assumption of independence typically found in studies such as \citep{schmiedel2024, boujemaa2022}, which model the outflow (demand) in transportation systems. The Gaussian copula enables a flexible correlation structure while maintaining Poisson marginals. Section~\ref{sec:nuisance} outlines the methodology for generating such data, including extensions for scenarios where researchers may need to incorporate autocorrelation or seasonal trends.

The performance of the SDDP algorithm, as presented in Table~\ref{tab:perf6x6}, shows a consistent improvement in solution quality with increasing iterations. Based on 100 problem instances, bias and CI ratio decrease from \(35.5 \%\) and \(21.5 \%\) at 500 iterations to \(15.7 \%\) and \(9.6 \%\) at 1,500 iterations, respectively, indicating greater accuracy in the value function approximation and reduced uncertainty in the in-sample objective value. However, this improvement comes with a substantial increase in computational cost, as the average solution time rises non-linearly from \(18.90\) seconds at 500 iterations to \(118.49\) seconds at 1,500 iterations. This trade-off between computational efficiency and solution quality highlights the diminishing returns of additional iterations, emphasizing the importance of selecting a balanced iteration count to achieve both accuracy and efficiency.

\begin{table}
    \centering
    \caption{SDDP performance measures on 100 instances of the problem.}
    \label{tab:perf6x6}
    \begin{tabular}{rrrrrr}
        \toprule
        Iterations & Solves & Bias $(\%)$ & CI ratio $(\%)$  & Avg time (sec) \\
        \midrule
         500  &  67,188 & \( 35.5 \pm 0.4 \) & \( 21.5 \pm 0.2 \) & \( 18.90 \pm 0.21 \) \\
        1,000 & 133,188 & \( 21.8 \pm 0.3 \) & \( 13.3 \pm 0.2 \) & \( 37.49 \pm 0.25 \) \\
        1,500 & 199,188 & \( 15.7 \pm 0.2 \) & \(  9.6 \pm 0.1 \) & \( 118.49 \pm 1.51 \) \\
        \bottomrule
    \end{tabular}
\end{table}

\begin{figure}
    \centering
    \includegraphics[width = \textwidth]{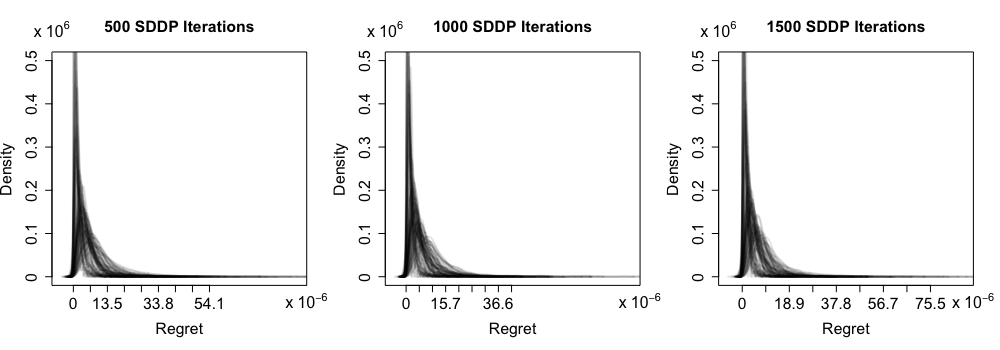}
    \caption{Distribution of regret recorded from 1,000 out-of-bag inflow/outflow scenarios simulated for 100 instances of the 12-stage stochastic transportation problem with 6 entry hubs, 6 exit hubs, and 20 carriers.}
    \label{fig:regret100x1000}
\end{figure}

The regret analysis in Figure~\ref{fig:regret100x1000} illustrates the distribution of regret across 1,000 out-of-bag inflow/outflow scenarios for 100 instances of the 12-stage stochastic transportation problem. Across all iteration counts—500, 1,000, and 1,500—the regret is consistently well-concentrated near zero, with minimal differences in the distributions. This suggests that increasing the number of iterations has little impact on reducing regret, as the algorithm performs robustly even at lower iteration counts. The marginal improvements observed at higher iterations are negligible, indicating that fewer iterations may suffice to achieve high-quality policies, offering a computationally efficient alternative without significant loss in performance.

The consistency in regret distributions across iteration counts highlights the diminishing returns of additional iterations in solving the stochastic transportation problem. The negligible improvements at higher iterations suggest that the algorithm effectively captures the critical structure of the value function early in training, achieving robust and near-optimal policies with fewer iterations. This raises important questions about the necessity of prolonged computation when high-quality solutions can be obtained efficiently. A more adaptive stopping criterion could balance performance metrics like regret against computational cost, ensuring scalable and resource-efficient solutions for large-scale stochastic optimization.

\section{Discussion}\label{sec:disco}

We investigated the MSTP within the context of drayage procurement for container logistics, addressing the challenges posed by uncertainty in cargo flows and market conditions. By leveraging SDDP, we demonstrated an efficient and scalable approach for deriving optimal operational policies across complex, realistic problem instances. The findings provide valuable insights into the influence of uncertainty, computational trade-offs, and the practical applicability of advanced stochastic optimization techniques in container logistics.


Our results emphasize the critical role of inflow uncertainty in driving total operational costs, as evidenced by the sensitivity analysis. Variations in cargo inflows significantly impact optimal allocation decisions, underscoring the importance of incorporating inflow-side uncertainty into optimization models. By contrast, spot rate variability exhibited a minimal influence on costs, suggesting that under current procurement strategies, reliance on spot markets remains limited. This result aligns with industry practice, where shippers prioritize long-term contracts to secure capacity and use the spot market as a buffer for flexibility. Therefore, the proposed approach supports procurement managers in better anticipating the operational impact of flow uncertainties while optimizing cost structures through effective contract and spot market utilization.

The analysis of volume allocation policies highlights the balance between contracted capacity and spot market flexibility. Contracted carriers play a dominant role in ensuring cost stability and meeting volume commitments, while opportunistic spot market procurement enhances operational efficiency. This insight provides a foundation for developing dynamic procurement strategies that can adapt to evolving uncertainties without over-relying on costly reactive measures.


The numerical experiments confirm the scalability and robustness of the SDDP approach in solving multistage stochastic problems of industrial scale. The algorithm achieves high-quality solutions with relatively low regret across out-of-bag scenarios, even at a moderate number of iterations. This indicates that SDDP effectively captures the value function structure early in the iterative process, reducing the need for prolonged computation. Importantly, the trade-off between computational cost and solution accuracy is made explicit through our analysis of bias, confidence intervals, and execution times.

The observed diminishing marginal improvements at higher iteration counts suggest opportunities to introduce adaptive stopping rules. Such rules could terminate iterations once a predefined accuracy threshold is achieved, balancing computational efficiency and solution fidelity. This finding is particularly relevant for large-scale systems, where computational resources are limited, and timely decision-making is essential.


Our work advances the field of truckload procurement by illustrating how the multivariate time series framework of \citet{Fokianos2020} can be employed to relax common independence assumptions across time stages or between flows (supply and demand). By leveraging copula-based techniques, we demonstrate a straightforward yet effective method for capturing both autocorrelations and cross-correlations in cargo dynamics, resulting in a more realistic representation of uncertainty. Additionally, we propose an alternative stopping criterion for SDDP iterations, offering a practical improvement that balances computational effort with solution quality in stochastic programming.



While the study offers valuable insights, further enhancements could strengthen its applicability. The computational experiments are based on specific parameter choices and assumptions about uncertainty distributions. Although the copula-based approach improves realism, validating the model with empirical data from container logistics operations would increase its robustness. Additionally, the current framework assumes fixed carrier capacities and stationary spot market rates, which could be extended to incorporate dynamic capacity decisions and pricing mechanisms. Future research could also consider operational constraints, such as routing decisions and time windows, to provide a more comprehensive drayage solution. Exploring alternative techniques, such as data-driven function approximations, may further enhance computational performance for large-scale problems.


In summary, this study demonstrates the effectiveness of SDDP in solving MSTP for truckload procurement under uncertainty. By addressing the computational challenges and highlighting the sensitivity of costs to flow uncertainty, we provide both theoretical contributions and practical insights for container logistics management. The findings pave the way for more robust and scalable procurement strategies, ensuring adaptability and cost efficiency in dynamic transportation networks.

\bibliographystyle{apalike}
\bibliography{refs.bib}

\end{document}